# Proving Time Bounds for Randomized Distributed Algorithms [*]
## (Extended Abstract)

Nancy Lynch[†]     Isaac Saias[‡]     Roberto Segala[§]

Laboratory for Computer Science
Massachusetts Institute of Technology
Cambridge, MA 02139

**Abstract**

A method of analyzing time bounds for randomized distributed algorithms is presented, in the context of a new and general framework for describing and reasoning about randomized algorithms. The method consists of proving auxiliary statements of the form $U \xrightarrow[p]{t} U'$, which means that whenever the algorithm begins in a state in set $U$, with probability $p$, it will reach a state in set $U'$ within time $t$. The power of the method is illustrated by its use in proving a constant upper bound on the expected time for some process to reach its critical region, in Lehmann and Rabin's Dining Philosophers algorithm.

**Contact Author:** Roberto Segala

[*]Supported by NSF grant CCR-89-15206, and CCR-92-25124, by DARPA contracts N00014-89-J-1988 and N00014-92-J-4033, and by ONR contract N00014-91-J-1046.
[†]e-mail: lynch@theory.mit.edu
[‡]e-mail: saias@theory.lcs.mit.edu
[§]e-mail: segala@theory.lcs.mit.edu

# 1 Introduction

Randomization has proved to be a useful tool in the design of distributed algorithms, sometimes yielding efficient solutions to problems that are inherently complex, or even unsolvable, in the setting of deterministic algorithms [?,?,?,?]. But this powerful tool has a price: even simple randomized algorithms can be extremely hard to verify and analyze. Because of this, many randomized distributed algorithms appear in the literature with only informal proofs of correctness, and only informal derivation of complexity bounds. In fact, it is sometimes hard for the reader to ascertain that the proofs and complexity bounds presented are really correct. Even where proofs are carefully and correctly done, the arguments tend to be ad hoc.

A key difficulty in reasoning about randomized algorithms is the fact that their executions usually contain a combination of nondeterministic and probabilistic choices, with subtle interactions between them. The probabilistic choices are typically only those that involve an explicit use of randomness by the algorithm (e.g., by using a random-number generator). All other choices (e.g., the order of process steps, the times at which requests arrive) are usually considered to be nondeterministic. It is customary to define an *adversary* as a way of modeling the entity that resolves the nondeterministic choices.[1] In defining an adversary, one must be especially careful to specify the *knowledge* of the execution that the adversary is permitted to use in resolving nondeterministic choices. This might range from no knowledge at all, in which case the adversary is said to be *oblivious*, to complete knowledge of the past execution (including past random choices).

Even after one has defined the desired notion of adversary, it is still not easy to carry out correctness proofs and complexity analyses for randomized algorithms; most existing proofs seem rather ad hoc. It would be useful to have a collection of general proof rules and methods, which could be established once and for all, and then applied in a reasonably systematic way to verify and analyze numerous algorithms. Some examples of work that has already been done on the development of such methods is [?,?,?]. The work of [?] presents a technique, based on progress functions defined on states, for establishing liveness properties for randomized algorithms; the work of [?,?] presents model checking techniques.

In this paper, we present such a new method: a way of proving *upper bounds on time* for randomized algorithms. Our method consists of proving auxiliary statements of the form $U \xrightarrow[p]{t} U'$, which means that whenever the algorithm begins in a state in set $U$, with probability $p$, it will reach a state in set $U'$ within time $t$. Of course, this method can only be used for randomized algorithms that include timing assumptions. A key theorem about our method is the composability of these $U \xrightarrow[p]{t} U'$ arrows, as expressed by Theorem 3.4. This composability result holds even in the case of (many classes of) non-oblivious adversaries.

We also present two complementary proof rules that help in reasoning about sets of distinct random choices. Independence arguments about such choices are often crucial to correctness proofs, yet there are subtle ways in which a non-oblivious adversary can introduce dependencies. For example, a non-oblivious adversary has the power to use the outcome of one random choice to decide whether to schedule another random choice. Our proof rules help to systematize certain kinds of reasoning about independence.

---

[1] In this paper, we ignore the possibility that the adversary itself uses randomness.

Our proof rules are presented in the context of a new and general formal framework [?] for describing and reasoning about randomized algorithms. This framework integrates randomness and nondeterminism into one model, and permits the modeling of timed as well as untimed systems. The model of [?] is, in turn, based on existing models for untimed and timed distributed systems [?, ?], and adopts many ideas from the probabilistic models of [?, ?].

In order to illustrate our method, we use it in this paper to prove an upper bound for Lehmann and Rabin's Dining Philosophers algorithm [?], in the face of an adversary with complete knowledge of the past. This upper bound asserts that $\mathcal{T} \xrightarrow[1/8]{13} \mathcal{C}$, where $\mathcal{T}$ is the set of states in which some process is in its trying region, while $\mathcal{C}$ is the set of states in which some process is in its critical region. That is, whenever the algorithm is in a state in which some process is in the trying region, with probability 1/8, within time 13, it will reach a state in which some process is in its critical region. This bound depends on the timing assumption that processes never wait more then time 1 between steps. A consequence of this claim is an upper bound (of 63) on the expected time for some process to reach its critical region.

For comparison, we note that [?] contains only proof sketches of the results claimed. The paper [?] contains a proof that Lehmann and Rabin's algorithm satisfies an *eventual* progress condition, in the presence of an adversary with complete knowledge of the past; this proof is carried out as an instance of Zuck and Pnueli's general method for proving liveness properties. Our results about this protocol can be regarded as a refinement of the results of Zuck and Pnueli, in that we obtain explicit constant time bounds rather than liveness properties.

The rest of the paper is organized as follows. Section 2 presents a simplified version of the model of [?]. Section 3 presents our main proof technique based on time-bound statements. Section 4 presents the additional proof rules for independence of distinct probabilistic choices. Section 5 presents the Lehmann-Rabin algorithm. Section 6.2 formalizes the algorithm in terms of the model of Section 2, and gives an overview of our time bound proof. Section 7 gives some concluding remarks. A separate appendix contains the details of the time bound proof.

## 2   The Model

In this section, we present the model that is used to formulate our proof technique. It is a simplified version of the probabilistic automaton model of [?]. Here we only give the parts of the model that we need to describe our proof method and its application to the Lehmann-Rabin algorithm; we refer the reader to the full version of this paper and to [?] for more details.

**Definition 2.1** A *probabilistic automaton*[2] $M$ consists of four components:

- a set *states*($M$) of states.
- a nonempty set *start*($M$) $\subseteq$ *states*($M$) of start states.

---

[2]In [?] the probabilistic automata of this definition are called *simple probabilistic automata*. This is because that paper also includes the case of randomized adversaries.

- an action signature $sig(M) = (ext(M), int(M))$ where $ext(M)$ and $int(M)$ are disjoint sets of external and internal actions, respectively.

- a transition relation $steps(M) \subseteq states(M) \times acts(M) \times Probs(states(M))$, where the set $Probs(states(M))$ is the set of probability spaces $(\Omega, \mathcal{F}, P)$ such that $\Omega \subseteq states(M)$ and $\mathcal{F} = 2^\Omega$. The last requirement is needed for technical convenience.

A probabilistic automaton is *fully probabilistic* if it has a unique start state and from each state there is at most one step enabled. ∎

Thus, a probabilistic automaton is a state machine with a labeled transition relation such that the state reached during a step is determined by some probability distribution. For example, the process of flipping a coin is represented by a step labeled with an action `flip` where the next state contains the outcome of the coin flip and is determined by a probability distribution over the two possible outcomes. A probabilistic automaton also allows nondeterministic choices over steps. An example of nondeterminism is the choice of which process takes the next step in a multi-process system.

An *execution fragment* $\alpha$ of a probabilistic automaton $M$ is a (finite or infinite) sequence of alternating states and actions starting with a state and, if the execution fragment is finite, ending in a state, $\alpha = s_0 a_1 s_1 a_2 s_2 \cdots$, where for each $i$ there exists a probability space $(\Omega, \mathcal{F}, P)$ such that $(s_i, a_{i+1}, (\Omega, \mathcal{F}, P)) \in steps(M)$ and $s_{i+1} \in \Omega$. Denote by $fstate(\alpha)$ the first state of $\alpha$ and, if $\alpha$ is finite, denote by $lstate(\alpha)$ the last state of $\alpha$. Furthermore, denote by $frag^*(M)$ and $frag(M)$ the sets of finite and all execution fragments of $M$, respectively. An *execution* is an execution fragment whose first state is a start state. Denote by $exec^*(M)$ and $exec(M)$ the sets of finite and all executions of $M$, respectively. A state $s$ of $M$ is *reachable* if there exists a finite execution of $M$ that ends in $s$. Denote by $rstates(M)$ the set of reachable states of $M$.

A finite execution fragment $\alpha_1 = s_0 a_1 s_1 \cdots a_n s_n$ of $M$ and an execution fragment $\alpha_2 = s_n a_{n+1} s_{n+1} \cdots$ of $M$ can be *concatenated*. In this case the concatenation, written $\alpha_1 \frown \alpha_2$, is the execution fragment $s_0 a_1 s_1 \cdots a_n s_n a_{n+1} s_{n+1} \cdots$. An execution fragment $\alpha_1$ of $M$ is a *prefix* of an execution fragment $\alpha_2$ of $M$, written $\alpha_1 \leq \alpha_2$, if either $\alpha_1 = \alpha_2$ or $\alpha_1$ is finite and there exists an execution fragment $\alpha_1'$ of $M$ such that $\alpha_2 = \alpha_1 \frown \alpha_1'$.

In order to study the probabilistic behavior of a probabilistic automaton, some mechanism to remove nondeterminism is necessary. To give an idea of why the nondeterministic behavior should be removed, consider a probabilistic automaton with three states $s_0, s_1, s_2$ and with two steps enabled from its start state $s_0$; the first step moves to state $s_1$ with probability $1/2$ and to $s_2$ with probability $1/2$; the second step moves to state $s_1$ with probability $1/3$ and to $s_2$ with probability $2/3$. What is the probability of reaching state $s_1$? The answer depends on how the nondeterminism between the two steps is resolved. If the first step is chosen, then the probability of reaching state $s_1$ is $1/2$; if the second step is chosen, then the probability of reaching state $s_1$ is $1/3$. We call the mechanism that removes the nondeterminism an *adversary*, because it is often viewed as trying to thwart the efforts of a system to reach its goals. In distributed systems the adversary is often called the *scheduler*, because its main job may be to decide which process should take the next step.

**Definition 2.2** An *adversary* for a probabilistic automaton $M$ is a function $\mathcal{A}$ taking a finite execution fragment of $M$ and giving back either nothing (represented as $\delta$) or one of the enabled steps of $M$ if there are any. Denote the set of adversaries for $M$ by $Advs_M$[3]. ∎

Once an adversary is chosen, a probabilistic automaton can run under the control of the chosen adversary. The result of the interaction is called an execution automaton. The definition of an execution automaton, given below, is rather complicated because an execution automaton must contain all the information about the different choices of the adversary, and thus the states of an execution automaton must contain the complete history of a probabilistic automaton. Note that there are no nondeterministic choices left in an execution automaton.

**Definition 2.3** An *execution automaton* $H$ of a probabilistic automaton $M$ is a fully probabilistic automaton such that

1. $states(H) \subseteq frag^*(M)$.

2. for each step $(\alpha, a, (\Omega, \mathcal{F}, P))$ of $H$ there is a step $(lstate(\alpha), a, (\Omega', \mathcal{F}', P'))$ of $M$, called the corresponding step, such that $\Omega = \{\alpha a s | s \in \Omega'\}$ and $P'[\alpha a s] = P[s]$ for each $s \in \Omega'$.

3. each state of $H$ is reachable, i.e., for each $\alpha \in states(H)$ there exists an execution of $H$ leading to state $\alpha$. ∎

**Definition 2.4** Given a probabilistic automaton $M$, an adversary $\mathcal{A} \in Advs_M$, and an execution fragment $\alpha \in frag^*(M)$, the execution $H(M, \mathcal{A}, \alpha)$ of $M$ under adversary $\mathcal{A}$ with starting fragment $\alpha$ is the execution automaton of $M$ whose start state is $\alpha$ and such that for each step $(\alpha', a, (\Omega, \mathcal{F}, P)) \in steps(H(M, \mathcal{A}, \alpha))$, its corresponding step is the step $\mathcal{A}(\alpha')$. ∎

Given an execution automaton $H$, an event is expressed by means of a set of maximal executions of $H$, where a maximal execution of $H$ is either infinite, or it is finite and its last state does not enable any step in $H$. For example, the event "eventually action $a$ occurs" is the set of maximal executions of $H$ where action $a$ does occur. A more formal definition follows. The sample space $\Omega_H$ is the set of maximal executions of $H$. The $\sigma$-algebra $\mathcal{F}_H$ is the smallest $\sigma$-algebra that contains the set of *rectangles* $R_\alpha$, consisting of the executions of $\Omega_H$ having $\alpha$ as a prefix[4]. The probability measure $P_H$ is the unique extension of the probability measure defined on rectangles as follows: $P_H[R_\alpha]$ is the product of the probabilities of each step of $H$ generating $\alpha$. In [?] it is shown that there is a unique probability measure having the property above, and thus $(\Omega_H, \mathcal{F}_H, P_H)$ is a well defined probability space. For the rest of this abstract we do not need to refer to this formal definition any more.

Events of $\mathcal{F}_H$ are not sufficient for the analysis of a probabilistic automaton. Events are defined over execution automata, but a probabilistic automaton may generate several execution automata depending on the adversary it interacts with. Thus a more general notion of event is needed that can deal with all execution automata. Specific examples are given in Section 3.

---

[3]In [?] the adversaries of this definition are denoted by $DAdvs_M$, where $D$ stands for *Deterministic*. The adversaries of [?] are allowed to use randomness.

[4]Note that a rectangle $R_\alpha$ can be used to express the fact that the finite execution $\alpha$ occurs.

**Definition 2.5** An *event schema* $e$ for a probabilistic automaton $M$ is a function associating an event of $\mathcal{F}_H$ with each execution automaton $H$ of $M$. ∎

We now discuss briefly a simple way to handle time within probabilistic automata. The idea is to add a time component to the states of a probabilistic automaton, to assume that the time at a start state is 0, to add a special non-visible action $\nu$ modeling the passage of time, and to add arbitrary time passage steps to each state. A time passage step should be non-probabilistic and should change only the time component of a state. This construction is called the *patient* construction in [?,?,?]. The reader interested in a more general extension to timed models is referred to [?].

We close this section with one final definition. Our time bound property for the Lehmann-Rabin algorithm states that if some process is in its trying region, then no matter how the steps of the system are scheduled, some process enters its critical region within time $t$ with probability at least $p$. However, this claim can only be valid if each process has sufficiently frequent chances to perform a step of its local program. Thus, we need a way to restrict the set of adversaries for a probabilistic automaton. The following definition provides a general way of doing this.

**Definition 2.6** An *adversary schema* for a probabilistic automaton $M$, denoted by $Advs$, is a subset of $Advs_M$. ∎

## 3 The Proof Method

In this section, we introduce our key statement $U \xrightarrow[p]{t}_{Advs} U'$ and the *composability theorem*, which is our main theorem about the proof method.

The meaning of the statement $U \xrightarrow[p]{t}_{Advs} U'$ is that, starting from any state of $U$ and under any adversary $\mathcal{A}$ of $Advs$, the probability of reaching a state of $U'$ within time $t$ is at least $p$. The suffix $Advs$ is omitted whenever we think it is clear from the context.

**Definition 3.1** Let $e_{U',t}$ be the event schema that, applied to an execution automaton $H$, returns the set of maximal executions $\alpha$ of $H$ where a state from $U'$ is reached in some state of $\alpha$ within time $t$. Then $U \xrightarrow[p]{t}_{Advs} U'$ iff for each $s \in U$ and each $\mathcal{A} \in Advs$, $P_{H(M,\mathcal{A},s)}[e_{U',t}(H(M,\mathcal{A},s))] \geq p$. ∎

**Proposition 3.2** *Let $U, U', U''$ be sets of states of a probabilistic automaton $M$. If $U \xrightarrow[p]{t} U'$, then $U \cup U'' \xrightarrow[p]{t} U' \cup U''$.* ∎

In order to compose time bound statements, we need a restriction for adversary schemas stating that the power of the adversary schema is not reduced if a prefix of the past history of the execution is not known. Most adversary schemas that appear in the literature satisfy this restriction.

**Definition 3.3** An adversary schema $Advs$ for a probabilistic automaton $M$ is *execution closed* if, for each $\mathcal{A} \in Advs$ and each finite execution fragment $\alpha \in \mathit{frag}^*(M)$, there exists an adversary $\mathcal{A}' \in Advs$ such that for each execution fragment $\alpha' \in \mathit{frag}^*(M)$ with $\mathit{lstate}(\alpha) = \mathit{fstate}(\alpha')$, $\mathcal{A}'(\alpha') = \mathcal{A}(\alpha \frown \alpha')$. ∎

**Theorem 3.4** *Let $Advs$ be an execution closed adversary schema for a probabilistic timed automaton $M$, and let $U, U', U''$ be sets of states of $M$.*
*If $U \xrightarrow[p_1]{t_1}_{Advs} U'$ and $U' \xrightarrow[p_2]{t_2}_{Advs} U''$, then $U \xrightarrow[p_1 p_2]{t_1+t_2}_{Advs} U''$.*

**Proof sketch.** Consider an adversary $\mathcal{A} \in Advs$ that acts on $M$ starting from a state $s$ of $U$. The execution automaton $H(M, \mathcal{A}, s)$ contains executions where a state from $U'$ is reached within time $t_1$. Consider one of those executions $\alpha$ and consider the part $H$ of $H(M, \mathcal{A}, s)$ after the first occurrence of a state from $U'$ in $\alpha$. The key idea of the proof is to use execution closure of $Advs$ to show that there is an adversary that generates $H$, to use $U' \xrightarrow[p_2]{t_2}_{Advs} U''$ to show that in $H$ a state from $U''$ is reached within time $t_2$ with probability at least $p_2$, and to integrate this last result in the computation of the probability of reaching a state from $U''$ in $H(M, \mathcal{A}, s)$ within time $t_1 + t_2$. ∎

## 4 Independence

**Example 4.1** Consider any distributed algorithm where each process is allowed to flip fair coins. It is common to say "If the next coin flip of process $P$ yields *head* and the next coin flip of process $Q$ yields *tail*, then some good property $\phi$ holds." Can we conclude that the probability for $\phi$ to hold is $1/4$? That is, can we assume that the coin flips of processes $P$ and $Q$ are independent? The two coin flips are indeed independent of each other, but the presence of non-oblivious adversaries may introduce some dependence. An adversary can schedule process $P$ to flip its coin and then schedule process $Q$ only if the coin flip of process $P$ yielded *head*. As a result, if both $P$ and $Q$ flip a coin, the probability that $P$ yields *head* and $Q$ yields *tail* is $1/2$. ∎

Thus, it is necessary to be extremely careful about independence assumptions. It is also important to pay attention to potential ambiguities of informal arguments. For example, does $\phi$ hold if process $P$ flips a coin yielding *head* and process $Q$ does not flip any coin? Certainly such an ambiguity can be avoided by expressing each event in a formal model.

In this section we present two event schemas that play a key role in the detailed time bound proof for the Lehmann-Rabin algorithm (cf. appendix), and we show some partial independence properties for them. The first event schema is a generalization of the informal statement of Example 4.1, where a coin flip is replaced by a generic action $a$, and where it is assumed that an event contains all the executions where $a$ is not scheduled; the second event schema is used to analyze the outcome of the first random draw that occurs among a fixed set of random draws. A consequence of the partial independence results that we show below is that under any adversary the property $\phi$ of Example 4.1 holds with probability at least $1/4$.

Let $(a, U)$ be a pair consisting of an action of $M$ and a set of states of $M$. The event schema FIRST$(a, U)$ is the function that, given an execution automaton $H$, returns the set of maximal executions of $H$ where either action $a$ does not occur, or action $a$ occurs and the state reached after the first occurrence of $a$ is a state of $U$. This event schema is used to express properties like "the $i^{\text{th}}$ coin yields left". For example $a$ can be flip and $U$ can be the set of states of $M$ where the result of the coin flip is left.

Let $(a_1, U_1), \ldots, (a_n, U_n)$ be a sequence of pairs consisting of an action of $M$ and a set of states of $M$ such that for each $i, j$, $1 \leq i < j \leq n$, $a_i \neq a_j$. Define the event schema NEXT$((a_1, U_1), \ldots, (a_n, U_n))$ to be the function that applied to an execution automaton $H$ gives the set of maximal executions of $H$ where either no action from $\{a_1, \ldots, a_n\}$ occurs, or at least one action from $\{a_1, \ldots, a_n\}$ occurs and, if $a_i$ is the first action that occurs, the state reached after the first occurrence of $a_i$ is in $U_i$. This kind of event schema is used to express properties like "the first coin that is flipped yields left."

**Proposition 4.2** Let $H$ be an execution automaton of a probabilistic automaton $M$. Furthermore, let $(a_1, U_1), \ldots, (a_n, U_n)$ be pairs consisting of an action of $M$ and a set of states of $M$ such that for each $i, j$, $1 \leq i < j \leq n$, $a_i \neq a_j$. Finally, let $p_1, \ldots, p_n$ be real numbers between 0 and 1 such that for each $i$, $1 \leq i \leq n$, and each step $(s, a, (\Omega, \mathcal{F}, P)) \in steps(M)$ with $a = a_i$, the probability $P[U_i \cap \Omega]$ is greater than or equal to $p_i$, i.e., $P[U_i \cap \Omega] \geq p_i$. Then

1. $P_H[(\text{FIRST}(a_1, U_1) \cap \cdots \cap \text{FIRST}(a_n, U_n))(H)] \geq p_1 \cdots p_n$,

2. $P_H[\text{NEXT}((a_1, U_1), \ldots, (a_n, U_n))(H)] \geq min(p_1, \ldots, p_n)$. ∎

## 5  The Lehmann-Rabin Algorithm

The Lehmann-Rabin algorithm is a randomized algorithm for the Dining Philosophers problem. This problem involves the allocation of $n$ resources among $n$ competing processes arranged in a ring. The resources are considered to be interspersed between the processes, and each process requires both its adjacent resources in order to reach its critical section. All processes are identical; the algorithm breaks symmetry by using randomization. The algorithm ensures the required exclusive possession of resources, and also ensures that, with probability 1, some process is always permitted to make progress into its critical region.

Figure 1 shows the code for a generic process $i$. The $n$ resources are represented by $n$ shared variables $\text{Res}_1, \ldots, \text{Res}_n$, each of which can assume values in \{free, taken\}. Each process $i$ ignores its own name, $i$, and the names, $\text{Res}_{i-1}$ and $\text{Res}_i$, of its adjacent resources. However, each process $i$ is able to refer to its adjacent resources by relative names: $\text{Res}_{(i,\text{left})}$ is the resource located to the left (clockwise), and $\text{Res}_{(i,\text{right})}$ is the resource to the right (counterclockwise) of $i$. Each process has a private variable $u_i$, which can assume a value in \{left, right\}, and is used to keep track of the first resource to be handled. For notational convenience we define an operator $opp$ that complements the value of its argument, i.e., $opp(\text{right}) = $ left and $opp(\text{left}) = $ right.

**Shared variables:** $\text{Res}_j \in \{\text{free}, \text{taken}\}$, $j = 1, \ldots, n$, initially free.

**Local variables:** $u_i \in \{\text{left}, \text{right}\}$, $i = 1, \ldots, n$

**Code for process $i$:**

| | | |
|---|---|---|
| 0. | try | ** beginning of Trying Section ** |
| 1. | $< u_i \leftarrow random >$ | ** choose left or right with equal probability ** |
| 2. | $<$ if $\text{Res}_{(i,u_i)} = $ free then | |
| | $\quad \text{Res}_{(i,u_i)} := $ taken | ** pick up first resource ** |
| | else goto 2. $>$ | |
| 3. | $<$ if $\text{Res}_{(i,opp(u_i))} = $ free then | |
| | $\quad \text{Res}_{(i,opp(u_i))} := $ taken; | ** pick up second resource ** |
| | $\quad$ goto 5. $>$ | |
| 4. | $< \text{Res}_{(i,u_i)} := $ free; goto 1.$>$ | ** put down first resource ** |
| 5. | crit | ** end of Trying Section ** |
| | ** Critical Section ** | |
| 6. | exit | ** beginning of Exit Section ** |
| 7. | $< u_i \leftarrow $ left or right | ** nondeterministic choice ** |
| | $\quad \text{Res}_{(i,opp(u_i))} := $ free $>$ | ** put down first resources ** |
| 8. | $< \text{Res}_{(i,u_i)} := $ free $>$ | ** put down second resources ** |
| 9. | rem | ** end of Exit Section ** |
| | ** Remainder Section ** | |

Figure 1: The Lehmann-Rabin algorithm

The atomic actions of the code are individual resource accesses, and they are represented in the form $<atomic\text{-}action>$ in Figure 1. We assume that at most one process has access to the shared resource at each time.

An informal description of the procedure is "choose a side randomly in each iteration. Wait for the resource on the chosen side, and, after getting it, just check *once* for the second resource. If this check succeeds, then proceed to the critical region. Otherwise, put down the first resource and try again with a new random choice."

Each process exchanges messages with an external user. In its idle state, a process is in its remainder region $R$. When triggered by a try message from the user, it enters the competition to get its resources: we say that it enters its trying region $T$. When the resources are obtained, it sends a crit message informing the user of the possession of these resources: we then say that the process is in its critical region $C$. When triggered by an exit message from the user, it begins relinquishing its resources: we then say that the process is in its exit region $E$. When the resources are relinquished its sends a rem message to the user and enters its remainder region.

# 6 Overview of the Proof

In this section, we give our high-level overview of the proof. We first introduce some notation, then sketch the proof strategy at a high level. Details of the proof appear in the Appendix.

## 6.1 Notation

In this section we define a probabilistic automaton $M$ which describes the system of Section 5. We assume that process $i + 1$ is on the right of process $i$ and that resource $\text{Res}_i$ is between processes $i$ and $i + 1$. We also identify labels modulo $n$ so that, for instance, process $n + 1$ coincides with process 1.

A state $s$ of $M$ is a tuple $(X_1, \ldots, X_n, \text{Res}_1, \ldots, \text{Res}_n, t)$ containing the local state $X_i$ of each process $i$, the value of each resource $\text{Res}_i$, and the current time $t$. Each local state $X_i$ is a pair $(pc_i, u_i)$ consisting of a program counter $pc_i$ and the local variable $u_i$. The program counter of each process keeps track of the current instruction in the code of Figure 1. Rather then representing the value of the program counter with a number, we use a more suggestive notation which is explained in the table below. Also, the execution of each instruction is represented by an action. Only actions $\text{try}_i$, $\text{crit}_i$, $\text{rem}_i$, $\text{exit}_i$ below are external actions.

| Number | $pc_i$ | Action name | Informal meaning |
|---|---|---|---|
| 0 | $R$ | $\text{try}_i$ | **R**eminder region |
| 1 | $F$ | $\text{flip}_i$ | Ready to **F**lip |
| 2 | $W$ | $\text{wait}_i$ | **W**aiting for first resource |
| 3 | $S$ | $\text{second}_i$ | Checking for **S**econd resource |
| 4 | $D$ | $\text{drop}_i$ | **D**ropping first resource |
| 5 | $P$ | $\text{crit}_i$ | **P**re-critical region |
| 6 | $C$ | $\text{exit}_i$ | **C**ritical region |
| 7 | $E_F$ | $\text{dropf}_i$ | **E**xit: drop **F**irst resource |
| 8 | $E_S$ | $\text{drops}_i$ | **E**xit: drop **S**econd resource |
| 9 | $E_R$ | $\text{rem}_i$ | **E**xit: move to **R**eminder region |

The start state of $M$ assigns the value $\texttt{free}$ to all the shared variables $\text{Res}_i$, the value $R$ to each program counter $pc_i$, and an arbitrary value to each variable $u_i$. The transition relation of $M$ is derived directly from Figure 1. For example, for each state where $pc_i = F$ there is an internal step $\text{flip}_i$ that changes $pc_i$ into $W$ and assigns $\texttt{left}$ to $u_i$ with probability $1/2$ and $\texttt{right}$ to $u_i$ with probability $1/2$; from each state where $X_i = (W, \texttt{left})$ there is a step $\text{wait}_i$ that does not change the state if $\text{Res}_{(i,\texttt{left})} = \texttt{taken}$, and changes $pc_i$ into $S$ and $\text{Res}_{(i,\texttt{left})}$ into $\texttt{taken}$ if $\text{Res}_{(i,\texttt{left})} = \texttt{free}$; for each state where $pc_i = E_F$ there are two steps with action $\text{dropf}_i$: one step sets $u_i$ to $\texttt{right}$ and makes $\text{Res}_{(i,\texttt{left})}$ free, and the other step sets $u_i$ to $\texttt{left}$ makes $\text{Res}_{(i,\texttt{right})}$ free. The two separate steps correspond to a nondeterministic choice that is left to the adversary. For time passage steps we assume that at any point an arbitrary amount of time can pass; thus, from each state of $M$ and each positive $\delta$ there is a time passage step that increases the time component of $\delta$ and does not affect the rest of the state.

The value of each pair $X_i$ can be represented concisely by the value of $pc_i$ and an arrow (to the left or to the right) which describes the value of $u_i$. Thus, informally, a process $i$ is in state $\underrightarrow{S}$ or $\underrightarrow{D}$ (resp. $\underleftarrow{S}$ or $\underleftarrow{D}$) when $i$ is in state $S$ or $D$ while holding its right (resp. left) resource; process $i$ is in state $\underrightarrow{W}$ (resp. $\underleftarrow{W}$) when $i$ is waiting for its right (resp. left) resource to become free; process $i$ is in state $\underrightarrow{E_S}$ (resp. $\underleftarrow{E_S}$) when $i$ is in its exit region and it is still holding its right (resp. left) resource. Sometimes we are interested in sets of pairs; for example, whenever $pc_i = F$ the value of $u_i$ is irrelevant. With the simple value of $pc_i$ we denote the set of the two pairs $\{(pc_i, \texttt{left}), (pc_i, \texttt{right})\}$. Finally, with the symbol # we denote any pair where $pc_i \in \{W, S, D\}$. The arrow notation is used as before.

For each state $s = (X_0, \ldots, X_{n-1}, \text{Res}_1, \ldots, \text{Res}_{n-1}, t)$ of $M$ we denote by $X_i(s)$ the pair $X_i$ and by $\text{Res}_i(s)$ the value of the shared variable $\text{Res}_i$ in state $s$. Also, for any set $S$ of states of a process $i$, we denote by $X_i \in S$, or alternatively $X_i = S$ the set of states $s$ of $M$ such that $X_i(s) \in S$. Sometimes we abuse notation in the sense that we write expressions like $X_i \in \{F, D\}$ with the meaning $X_i \in F \cup D$. Finally, we write $X_i = E$ for $X_i = \{E_F, E_S, E_R\}$, and we write $X_i = T$ for $X_i \in \{F, W, S, D, P\}$.

A first basic lemma states that a reachable state of $M$ is uniquely determined by the local states its processes and the current time. Based on this lemma, our further specifications of state sets will not refer to the shared variables; however, we consider only reachable states for the analysis. The proof of the lemma is a standard proof of invariants.

**Lemma 6.1** *For each reachable state $s$ of $M$ and each $i$, $1 \leq i \leq n$, $\text{Res}_i = \texttt{taken}$ iff $X_i(s) \in \{\underrightarrow{S}, \underrightarrow{D}, P, C, E_F, \underrightarrow{E_S}\}$ or $X_{i+1}(s) \in \{\underleftarrow{S}, \underleftarrow{D}, P, C, E_F, \underleftarrow{E_S}\}$. Moreover, for each reachable state $s$ of $M$ and each $i$, $1 \leq i \leq n$, it is not the case that $X_i(s) \in \{\underrightarrow{S}, \underrightarrow{D}, P, C, E_F, \underrightarrow{E_S}\}$ and $X_{i+1}(s) \in \{\underleftarrow{S}, \underleftarrow{D}, P, C, E_F, \underleftarrow{E_S}\}$, i.e., only one process at a time can hold one resource.* ∎

## 6.2 Proof Sketch

In this section we show that the RL-algorithm guarantees time bounded progress, i.e., that from every state where some process is in its trying region, some process subsequently enters its critical region within an expected constant time bound. We assume that each process that is ready to perform a step does so within time 1: process $i$ is ready to perform a step whenever it enables an action different from $\texttt{try}_i$ or $\texttt{exit}_i$. Actions $\texttt{try}_i$ and $\texttt{exit}_i$ are supposed to be under the control of the user, and hence, by assumption, under the control of the adversary.

Formally, consider the probabilistic timed automaton $M$ of Section 6.1. Define $Unit - Time$ to be the set of adversaries $\mathcal{A}$ for $M$ having the properties that, for every finite execution fragment $\alpha$ of $M$ and every execution $\alpha'$ of $H(M, \mathcal{A}, \alpha)$, 1) the time in $\alpha'$ is not bounded and 2) for every process $i$ and every state of $\alpha'$ enabling an action of process $i$ different from $\texttt{try}_i$ and $\texttt{exit}_i$, there exists a step in $\alpha'$ involving process $i$ within time 1. Then $Unit - Time$ is execution-closed according to Definition 3.3. An informal justification of this fact is that the constraint that each ready process is scheduled within time 1 knowing that $\alpha \frown \alpha'$ has occurred only reinforces the constraint that each ready process is scheduled within time 1 knowing that $\alpha'$ has occurred. Let

$$\mathcal{T} \triangleq \{s \in rstates(M) \mid \exists_i X_i(s) \in \{T\}\}$$

denote the sets of reachable states of $M$ where some process is in its trying region, and let

$$\mathcal{C} \triangleq \{s \in rstates(M) \mid \exists_i X_i(s) = C\}$$

denote the sets of reachable states of $M$ where some process is in its critical region. We show that

$$\mathcal{T} \xrightarrow[1/8]{13}_{Unit-Time} \mathcal{C},$$

i.e., that, starting from any reachable state where some process is in its trying region, for all the adversaries of $Unit - Time$, with probability at least $1/8$, some process enters its critical region within time $13$. Note that this property is trivially satisfied if some process is initially in its critical region.

Our proof is divided into several phases, each one concerned with the property of making a partial time bounded progress toward a "success state", i.e., a state of $\mathcal{C}$. The sets of states associated with the different phases are expressed in terms of $\mathcal{T}, \mathcal{RT}, \mathcal{F}, \mathcal{G}, \mathcal{P}$, and $\mathcal{C}$. Here,

$$\mathcal{RT} \triangleq \{s \in \mathcal{T} \mid \forall_i X_i(s) \in \{E_R, R, T\}\}$$

is the set of states where at least one process is in its trying region and where no process is in its critical region or holds resources while being in its exit region.

$$\mathcal{F} \triangleq \{s \in \mathcal{RT} \mid \exists_i X_i(s) = F\}$$

is the set of states of $\mathcal{RT}$ where some process is ready to flip a coin.

$$\mathcal{P} \triangleq \{s \in rstates(M) \mid \exists_i X_i(s) = P\}$$

is the sets of reachable states of $M$ where some process is in its pre-critical region. The set $\mathcal{G}$ is the most important for the analysis. It parallels the set of "Good Pairs" in [?] or the set described in Lemma 4 of [?]. To motivate the definition, we define the following notions. We say that a process $i$ is *committed* if $X_i \in \{W, S\}$, and that a process $i$ *potentially controls* $Res_i$ (resp. $Res_{i-1}$) if $X_i \in \{\overrightarrow{W}, \overrightarrow{S}, \overrightarrow{D}\}$ (resp. $X_i \in \{\overleftarrow{W}, \overleftarrow{S}, \overleftarrow{D}\}$). Informally said, a state in $\mathcal{RT}$ is in $\mathcal{G}$ if and only if there is a committed process whose second resource is not potentially controlled by another process. Such a process is called a *good* process. Formally,

$$\mathcal{G} \triangleq \{s \in \mathcal{RT} \mid \exists_i \ X_i(s) \in \{\overleftarrow{W}, \overleftarrow{S}\} \text{ and } X_{i+1}(s) \in \{E_R, R, F, \overleftarrow{\#}\}, \text{ or}$$
$$X_i(s) \in \{\overrightarrow{W}, \overrightarrow{S}\} \text{ and } X_{i-1}(s) \in \{E_R, R, F, \overrightarrow{\#}\}\}$$

Reaching a state of $\mathcal{G}$ is a substantial progress toward reaching a state of $\mathcal{C}$. Actually, the proof of Proposition A.11 establishes that, if $i$ a is good process, then, with probability $1/4$, one of the three processes $i - 1, i$ and $i + 1$ soon succeeds in getting its two resources. The hard part is to establish that, with constant probability, within a *constant* time, $\mathcal{G}$ is reached from any state in $\mathcal{T}$. A close inspection of the proof given in [?] shows that, there, the timed version of the techniques used is unable to deliver this result. The phases of our proof are formally described below.

$$\mathcal{T} \xrightarrow{2} \mathcal{RT} \cup \mathcal{C} \qquad \text{(Proposition A.3)},$$
$$\mathcal{RT} \xrightarrow{3} \mathcal{F} \cup \mathcal{G} \cup \mathcal{P} \qquad \text{(Proposition A.15)},$$
$$\mathcal{F} \xrightarrow[1/2]{2} \mathcal{G} \cup \mathcal{P} \qquad \text{(Proposition A.14)},$$
$$\mathcal{G} \xrightarrow[1/4]{5} \mathcal{P} \qquad \text{(Proposition A.11)},$$
$$\mathcal{P} \xrightarrow{1} \mathcal{C} \qquad \text{(Proposition A.1)}.$$

The first statement states that, within time 2, every process in its exit region relinquishes its resources. By combining the statements above by means of Proposition 3.2 and Theorem 3.4 we obtain

$$\mathcal{T} \xrightarrow[1/8]{13} \mathcal{C},$$

which is the property that was to be proven. Using the results of the proof summary above, we can furthermore derive a constant upper bound on the expected time required to reach a state of $\mathcal{C}$ when departing from a state of $\mathcal{T}$. Note that, departing from a state in $\mathcal{RT}$, with probability at least $1/8$, $\mathcal{P}$ is reached in time (at most) 10; with probability at most $1/2$, time 5 is spent before failing to reach $\mathcal{G} \cup \mathcal{P}$ ("failure at the third arrow"); with probability at most $7/8$, time 10 is spent before failing to reach $\mathcal{P}$ ("failure at the fourth arrow"). If failure occurs, then the state is back into $\mathcal{RT}$. Let $V$ denote a random variable satisfying the following induction

$$V = 1/8 \cdot 10 + 1/2\,(5 + V_1) + 3/8\,(10 + V_2)\,,$$

where $V_1$ and $V_2$ are random variables having the same distribution as $V$. The previous discussion shows that the expected time spent from $\mathcal{RT}$ to $\mathcal{P}$ is at most $E[V]$. By taking expectation in the previous equation, and using that $E[V] = E[V_1] = E[V_2]$, we obtain that $E[V] = 60$ is an upper bound on the expected time spent from $\mathcal{RT}$ to $\mathcal{P}$, and that, consequently, the expected time for progress starting from a state of $\mathcal{T}$ is at most 63.

## 7 Concluding Remarks

This paper has presented a formal model and a formal proof technique for the estimation of time performance of randomized algorithms that run under the control of general classes of adversaries. The salient aspect of this technique is to prove probabilistic time bounded progress properties and to compose them by means of a powerful composability theorem. The power of the proof method has been illustrated by proving a constant upper bound on the expected time for progress in the Lehmann-Rabin Dining Philosophers algorithm.

We believe that this technique is applicable towards the time analysis of many randomized protocols. It is desirable that the general model and this technique be used for the analysis of other algorithms, so that the power of the method can be tested and/or increased by means of other additional tools. In particular, it is very likely that new event schemas and partial independence results similar to those of Section 4 can be developed.

The specific results about the Lehmann-Rabin Dining Philosophers algorithm can be complemented and extended in many ways. We cite two. First, it would be very satisfying to derive a non trivial lower bound on the time for progress, which should be lower than our upper bound since the upper bound could be easily improved by means of a finer analysis. Second, it would be interesting to consider topologies that are more general than rings.

# Appendix

## A  The Detailed Proof

In this appendix we prove the five relations used in Section 6.2. However, for the sake of clarity, we do not prove the relations in the order they were presented. Throughout the proof we abuse notation by writing events of the kind $\text{FIRST}(\texttt{flip}_i, \texttt{left})$ meaning the event schema $\text{FIRST}(\texttt{flip}_i, \{s \in states(M) \mid X_i(s) = \underleftarrow{W}\})$.

**Proposition A.1** *If some process is in $P$, then, within time 1, it enters $C$, i.e.,*

$$\mathcal{P} \xrightarrow{1}_{1} \mathcal{C}.$$

**Proof.** This step corresponds to the action $\texttt{crit}$: within time 1, process $i$ informs the user that the critical region is free. ∎

**Lemma A.2** *If some process is in its Exit region then, within time 3, it will enter $R$.*

**Proof.** The process needs to take first two steps to relinquish its two resources, and then one step to send a $\texttt{rem}$ message to the user. ∎

**Proposition A.3** $\mathcal{T} \xrightarrow{2} \mathcal{RT} \cup \mathcal{C}$.

**Proof.** From Lemma A.2 within time 2 every process that begins in $E_F$ or $E_S$ relinquishes its resources. If no process begins in $C$ or enters $C$ in the meantime, then the state reached at this point is a state of $\mathcal{RT}$; otherwise, the starting state or the state reached when the first process enters $C$ is a state of $\mathcal{C}$. ∎

We now turn to the proof of $\mathcal{G} \xrightarrow{5}_{1/4} \mathcal{P}$. The following lemmas form a detailed cases analysis of the different situations that can arise in states of $\mathcal{G}$. Informally, each lemma shows that some event of the form of Proposition 4.2 is a sub-event of the properties of reaching some other state.

**Lemma A.4**

1. *Assume that $X_{i-1} \in \{E_R, R, F\}$ and $X_i = \underleftarrow{W}$. If $\text{FIRST}(\texttt{flip}_{i-1}, \texttt{left})$, then, within time 1, either $X_{i-1} = P$ or $X_i = S$.*

2. *Assume that $X_{i-1} = D$ and $X_i = \underleftarrow{W}$. If $\text{FIRST}(\texttt{flip}_{i-1}, \texttt{left})$, then, within time 2, either $X_{i-1} = P$ or $X_i = S$.*

3. *Assume that $X_{i-1} = S$ and $X_i = \underleftarrow{W}$. If $\text{FIRST}(\texttt{flip}_{i-1}, \texttt{left})$, then, within time 3, either $X_{i-1} = P$ or $X_i = S$.*

4. *Assume that $X_{i-1} = W$ and $X_i = \underleftarrow{W}$. If* FIRST($\texttt{flip}_{i-1}$, $\texttt{left}$), *then, within time 4, either $X_{i-1} = P$ or $X_i = S$.*

**Proof.** The four proofs start in the same way. Let $s$ be a state of $M$ satisfying the respective properties of items *1* or *2* or *3* or *4*. Let $f$ be an adversary of $Unit - Time$, and let $\alpha$ be the execution of $M$ that corresponds to an execution of $H(M, \{s\}, f)$ where the result of the first coin flip of process $i - 1$ is $\texttt{left}$.

1. By hypothesis, $i - 1$ does not hold any resource at the beginning of $\alpha$ and has to obtain $\text{Res}_{i-2}$ (its left resource) before pursuing $\text{Res}_{i-1}$. Within time 1, $i$ takes a step in $\alpha$. If $i - 1$ does not hold $\text{Res}_{i-1}$ when $i$ takes this step, then $i$ progresses into configuration $S$. If not, it must be the case that $i - 1$ succeeded in getting it in the meanwhile. But, in this case, $\text{Res}_{i-1}$ was the second resource needed by $i - 1$ and $i - 1$ therefore entered $P$.

2. If $X_i = S$ within time 1, then we are done. Otherwise, after one unit of time, $X_i$ is still equal to $\underleftarrow{W}$, i.e., $X_i(s') = \underleftarrow{W}$ for all states $s'$ reached in time 1. However, also process $i - 1$ takes a step within time 1. Let $\alpha = \alpha_1 \frown \alpha_2$ such that the last step of $\alpha_1$ is the first step taken by process $i - 1$. Then $X_{i-1}(\textit{fstate}(\alpha_2)) = F$ and $X_i(\textit{fstate}(\alpha_2)) = \underleftarrow{W}$. Since process $i - 1$ did not flip any coin during $\alpha_1$, from the execution closure of $Unit - Time$ and item *1* we conclude.

3. If $X_i = S$ within time 1, then we are done. Otherwise, after one unit of time, $X_i$ is still equal to $\underleftarrow{W}$, i.e., $X_i(s') = \underleftarrow{W}$ for all states $s'$ reached in time 1. However, also process $i - 1$ takes a step within time 1. Let $\alpha = \alpha_1 \frown \alpha_2$ such that the last step of $\alpha_1$ is the first step taken by process $i - 1$. If $X_{i-1}(\textit{fstate}(\alpha_2)) = P$ then we are also done. Otherwise it must be the case that $X_{i-1}(\textit{fstate}(\alpha_2)) = D$ and $X_i(\textit{fstate}(\alpha_2)) = \underleftarrow{W}$. Since process $i - 1$ did not flip any coin during $\alpha_1$, from the execution closure of $Unit - Time$ and item *2* we conclude.

4. If $X_i = S$ within time 1, then we are done. Otherwise, after one unit of time, $X_i$ is still equal to $\underleftarrow{W}$, i.e., $X_i(s') = \underleftarrow{W}$ for all states $s'$ reached in time 1. However, since within time 1 process $i$ checks its left resource and fails, process $i - 1$ gets its right resource within time 1, and hence reaches at least state $S$. Let $\alpha = \alpha_1 \frown \alpha_2$ where the last step of $\alpha_1$ is the first step of $\alpha$ leading process $i - 1$ to state $S$. Then $X_{i-1}(\textit{fstate}(\alpha_2)) = S$ and $X_i(\textit{fstate}(\alpha_2)) = \underleftarrow{W}$. Since process $i - 1$ did not flip any coin during $\alpha_1$, from the execution closure of $Unit - Time$ and item *3* we conclude. ∎

**Lemma A.5** *Assume that $X_{i-1} \in \{E_R, R, T\}$ and $X_i = \underleftarrow{W}$. If* FIRST($\texttt{flip}_{i-1}$, $\texttt{left}$), *then, within time 4, either $X_{i-1} = P$ or $X_i = S$.*

**Proof.** The lemma follows immediately from Lemma A.4 after observing that $X_{i-1} \in \{E_R, R, T\}$ means $X_{i-1} \in \{E_R, R, F, W, S, D, P\}$. ∎

The next lemma is a useful tool for the proofs of Lemmas A.7, A.8, and A.9.

**Lemma A.6** *Assume that $X_i \in \{\overleftrightarrow{W}, \overleftarrow{S}\}$ or $X_i \in \{E_R, R, F, \underrightarrow{D}\}$ with $\text{FIRST}(\text{flip}_i, \text{left})$, and assume that $X_{i+1} \in \{\overleftrightarrow{W}, \overrightarrow{S}\}$ or $X_{i+1} \in \{E_R, R, F, \underrightarrow{D}\}$ with $\text{FIRST}(\text{flip}_{i+1}, \text{right})$. Then the first of the two processes $i$ or $i+1$ testing its second resource enters $P$ after having performed this test (if this time ever comes).*

**Proof.** By Lemma 6.1 $\text{Res}_i$ is free. Moreover, $\text{Res}_i$ is the second resource needed by both $i$ and $i+1$. Whichever tests for it first gets it and enters $P$. ∎

**Lemma A.7** *If $X_i = \underrightarrow{S}$ and $X_{i+1} \in \{\overleftrightarrow{W}, \overrightarrow{S}\}$ then, within time 1, one of the two processes $i$ or $i+1$ enters $P$. The same result holds if $X_i \in \{\overleftrightarrow{W}, \overleftarrow{S}\}$ and $X_{i+1} = \underrightarrow{S}$.*

**Proof.** Being in state $S$, process $i$ tests its second resource within time 1. An application of Lemma A.6 finishes the proof. ∎

**Lemma A.8** *Assume that $X_i = \underrightarrow{S}$ and $X_{i+1} \in \{E_R, R, F, \underrightarrow{D}\}$. If $\text{FIRST}(\text{flip}_{i+1}, \text{right})$, then, within time 1, one of the two processes $i$ or $i+1$ enters $P$. The same result holds if $X_i \in \{E_R, R, F, D\}$, $X_{i+1} = \underrightarrow{S}$ and $\text{FIRST}(\text{flip}_i, \text{left})$.*

**Proof.** Being in state $S$, process $i$ tests its second resource within time 1. An application of Lemma A.6 finishes the proof. ∎

**Lemma A.9** *Assume that $X_{i-1} \in \{E_R, R, T\}$, $X_i = \overleftarrow{W}$, and $X_{i+1} \in \{E_R, R, F, \underrightarrow{W}, \underrightarrow{D}\}$. If $\text{FIRST}(\text{flip}_{i-1}, \text{left})$ and $\text{FIRST}(\text{flip}_{i+1}, \text{right})$, then within time 5 one of the three processes $i-1$, $i$ or $i+1$ enters $P$.*

**Proof.** Let $s$ be a state of $M$ such that $X_{i-1}(s) \in \{E_R, R, T\}$, $X_i(s) = \overleftarrow{W}$, and $X_{i+1}(s) \in \{E_R, R, F, \underrightarrow{W}, \underrightarrow{D}\}$. Let $f$ be an adversary of $Unit - Time$, and let $\alpha$ be the execution of $M$ that corresponds to an execution of $H(M, \{s\}, f)$ where the result of the first coin flip of process $i-1$ is $\text{left}$ and the result of the first coin flip of process $i+1$ is $\text{right}$. By Lemma A.5, within time 4 either process $i-1$ reaches configuration $P$ in $\alpha$ or process $i$ reaches configuration $\underleftarrow{S}$ in $\alpha$. If $i-1$ reaches configuration $P$, then we are done. If not, then let $\alpha = \alpha_1 \frown \alpha_2$ such that $\text{lstate}(\alpha_1)$ is the first state $s'$ of $\alpha$ with $X_i(s') = \underleftarrow{S}$. If $i+1$ enters $P$ before the end of $\alpha_1$, then we are done. Otherwise, $X_{i+1}(\text{fstate}(\alpha_2))$ is either in $\{\underrightarrow{W}, \underrightarrow{S}\}$ or it is in $\{E_R, R, F, \underrightarrow{D}\}$ and process $i+1$ has not flipped any coin yet in $\alpha$. From execution closure of $Unit - Time$ we can then apply Lemma A.6: within one more time process $i$ tests its second resource and by Lemma A.6 process $i$ enters $P$ if process $i+1$ did not check its second resource in the meantime. If process $i+1$ checks its second resource before process $i$ does the same, then by Lemma A.6 process $i+1$ enters $P$. Since process $i$ checks its second resource within time 1, process $i+1$ enters $P$ within time 1. ∎

**Lemma A.10** *Assume that $X_i \in \{E_R, R, F, \underleftarrow{W}, \underrightarrow{D}\}$, $X_{i+1} = \underrightarrow{W}$, and $X_{i+2} \in \{E_R, R, T\}$. If $\text{FIRST}(\text{flip}_i, \text{left})$ and $\text{FIRST}(\text{flip}_{i+2}, \text{right})$, then within time 5 one of the three processes $i$, $i+1$ or $i+2$, enters $P$.*

**Proof.** The proof is analogous to the one of Lemma A.9. This lemma is essentially the symmetric case of Lemma A.9. ∎

**Proposition A.11** *Starting from a global configuration in $\mathcal{G}$, then, with probability at least $1/4$ and within time at most $5$, some process enters $P$. Equivalently:*

$$\mathcal{G} \xrightarrow[1/4]{5} \mathcal{P}.$$

**Proof.** Lemmas A.7 and A.8 jointly treat the case where $X_i = \underleftarrow{S}$ and $X_{i+1} \in \{E_R, R, F, \underrightarrow{\#}\}$ and the symmetric case where $X_i \in \{E_R, R, F, \underleftarrow{\#}\}$ and $X_{i+1} = \underrightarrow{S}$; Lemmas A.9 and A.10 jointly treat the case where $X_i = \underleftarrow{W}$ and $X_{i+1} \in \{E_R, R, F, \underrightarrow{W}, \underrightarrow{D}\}$ and the symmetric case where $X_i \in \{E_R, R, F, \underleftarrow{W}, \underleftarrow{D}\}$ and $X_{i+1} = \underrightarrow{W}$.

Specifically, each lemma shows that a compound event of the kind FIRST($\texttt{flip}_i, x$) and FIRST($\texttt{flip}_j, y$) leads to $\mathcal{P}$. Each of the basic events FIRST($\texttt{flip}_i, x$) has probability $1/2$. From Proposition 4.2 each of the compound events has probability at least $1/4$. Thus the probability of reaching $\mathcal{P}$ within time $5$ is at least $1/4$. ∎

We now turn to $\mathcal{F} \xrightarrow[1/2]{2} \mathcal{G} \cup \mathcal{P}$. The proof is divided in two parts and constitute the global argument of the proof of progress.

**Lemma A.12** *Start with a state $s$ of $\mathcal{F}$. If there exists a process $i$ for which $X_i(s) = F$ and $(X_{i-1}, X_{i+1}) \neq (\underrightarrow{\#}, \underleftarrow{\#})$, then, with probability at least $1/2$ a state of $\mathcal{G} \cup \mathcal{P}$ is reached within time $1$.*

**Proof.** If $s \in \mathcal{G} \cup \mathcal{P}$, then the result is trivial. Let $s$ be a state of $\mathcal{F} - (\mathcal{G} \cup \mathcal{P})$ and let $i$ be such that $X_i(s) = F$ and $(X_{i-1}, X_{i+1}) \neq (\underrightarrow{\#}, \underleftarrow{\#})$. Assume without loss of generality that $X_{i+1} \neq \underleftarrow{\#}$, i.e., $X_{i+1} \in \{E_R, R, F, \underrightarrow{\#}\}$. The case for $X_{i-1} \neq \underrightarrow{\#}$ is similar. Furthermore, we can assume that $X_{i+1} \in \{E_R, R, F, \underrightarrow{D}\}$ since if $X_{i+1} \in \{\underrightarrow{W}, \underrightarrow{S}\}$ then $s$ is already in $\mathcal{G}$.

We show that the event NEXT(($\texttt{flip}_i, \texttt{left}$), ($\texttt{flip}_{i+1}, \texttt{right}$)), which by Proposition 4.2 has probability at least $1/2$, leads in time at most $1$ to a state of $\mathcal{G} \cup \mathcal{P}$. Let $f$ be an adversary of $Unit - Time$, and let $\alpha$ be the execution of $M$ that corresponds to an execution of $H(M, \{s\}, f)$ where if process $i$ flips before process $i+1$ then process $i$ flips left, and if process $i+1$ flips before process $i$ then process $i+1$ flips right.

Within time $1$, $i$ takes one step and reaches $W$. Let $j \in \{i, i+1\}$ be the first of $i$ and $i+1$ that reaches $W$ and let $s_1$ be the state reached after the first time process $j$ reaches $W$. If some process reached $P$ in the meantime, then we are done. Otherwise there are two cases to consider. If $j = i$, then, $\texttt{flip}_i$ gives $\texttt{left}$ and $X_i(s_1) = \underleftarrow{W}$ whereas $X_{i+1}$ is (still) in $\{E_R, R, F, \underrightarrow{D}\}$. Therefore, $s_1 \in \mathcal{G}$. If $j = i+1$, then $\texttt{flip}_{i+1}$ gives $\texttt{right}$ and $X_{i+1}(s_1) = \underrightarrow{W}$ whereas $X_i(s_1)$ is (still) $F$. Therefore, $s_1 \in \mathcal{G}$. ∎

**Lemma A.13** *Start with a state $s$ of $\mathcal{F}$. Assume that there exists a process $i$ for which $X_i(s) = F$ and for which $(X_{i-1}(s), X_{i+1}(s)) = (\underrightarrow{\#}, \underleftarrow{\#})$. Then, with probability at least $1/2$, within time $2$, a state of $\mathcal{G} \cup \mathcal{P}$ is reached.*

**Proof.** The hypothesis can be summarized into the form $(X_{i-1}(s), X_i(s), X_{i+1}(s)) = (\underset{\rightarrow}{\#}, F, \underset{\leftarrow}{\#})$. Since $i-1$ and $i+1$ point in different directions, by moving to the right of $i+1$ there is a process $k$ pointing to the left such that process $k+1$ either points to the right or is in $\{E_R, R, F, P\}$, i.e., $X_k(s) \in \{\underset{\leftarrow}{W}, \underset{\leftarrow}{S}, \underset{\leftarrow}{D}\}$ and $X_{k+1}(s) \in \{E_R, R, F, \underset{\rightarrow}{W}, \underset{\rightarrow}{S}, \underset{\rightarrow}{D}, P\}$. If $X_k(s) \in \{\underset{\leftarrow}{W}, \underset{\leftarrow}{S}\}$ and $X_{k+1}(s) \neq P$ then $s \in \mathcal{G}$ and we are done; if $X_{k+1}(s) = P$ then $s \in \mathcal{P}$ and we are done. Thus, we can restrict our attention to the case where $X_k(s) = \underset{\leftarrow}{D}$.

We show that the event NEXT$((\texttt{flip}_k, \texttt{left}), (\texttt{flip}_{k+1}, \texttt{right}))$, which by Proposition 4.2 has probability at least $1/2$, leads in time at most 2 to $\mathcal{G} \cup \mathcal{P}$. Let $f$ be an adversary of $Unit - Time$, and let $\alpha$ be the execution of $M$ that corresponds to an execution of $H(M, \{s\}, f)$ where if process $k$ flips before process $k+1$ then process $k$ flips left, and if process $k+1$ flips before process $k$ then process $k+1$ flips right.

Within time 2, process $k$ takes at least two steps and hence goes to configuration $W$. Let $j \in \{k, k+1\}$ be the first of $k$ and $k+1$ that reaches $W$ and let $s_1$ be the state reached after the first time process $j$ reaches $W$. If some process reached $P$ in the meantime, then we are done. Otherwise there are two cases to consider. If $j = k$, then, $\texttt{flip}_k$ gives $\texttt{left}$ and $X_k(s_1) = \underset{\leftarrow}{W}$ whereas $X_{k+1}$ is (still) in $\{E_R, R, F, \underset{\rightarrow}{\#}\}$. Therefore, $s_1 \in \mathcal{G}$. If $j = k+1$, then $\texttt{flip}_{k+1}$ gives $\texttt{right}$ and $X_{k+1}(s_1) = \underset{\rightarrow}{W}$ whereas $X_k(s_1)$ is (still) in $\{\underset{\leftarrow}{D}, F\}$. Therefore, $s_1 \in \mathcal{G}$. ∎

**Proposition A.14** *Start with a state $s$ of $\mathcal{F}$. Then, with probability at least $1/2$, within time 2, a state of $\mathcal{G} \cup \mathcal{P}$ is reached. Equivalently:*

$$\mathcal{F} \xrightarrow[1/2]{2} \mathcal{G} \cup \mathcal{P}.$$

**Proof.** The hypothesis of Lemmas A.12 and A.13 form a partition of $\mathcal{F}$. ∎

Finally, we prove $\mathcal{RT} \xrightarrow{3} \mathcal{F} \cup \mathcal{G} \cup \mathcal{P}$.

**Proposition A.15** *Starting from a state $s$ of $\mathcal{RT}$, then, within time 3, a state of $\mathcal{F} \cup \mathcal{G} \cup \mathcal{P}$ is reached. Equivalently:*

$$\mathcal{RT} \xrightarrow{3} \mathcal{F} \cup \mathcal{G} \cup \mathcal{P}.$$

**Proof.** Let $s$ be a state of $\mathcal{RT}$. If $s \in \mathcal{F} \cup \mathcal{G} \cup \mathcal{P}$, then we are trivially done. Suppose that $s \notin \mathcal{F} \cup \mathcal{G} \cup \mathcal{P}$. Then in $s$ each process is in $\{E_R, R, W, S, D\}$ and there exists at least process in $\{W, S, D\}$. Let $f$ be an adversary of $Unit - Time$, and let $\alpha$ be the execution of $M$ that corresponds to an execution of $H(M, \{s\}, f)$.

We first argue that within time 1 some process reaches a state of $\{S, D, F\}$ in $\alpha$. This is trivially true if in state $s$ there is some process in $\{S, D\}$. If this is not the case, then all processes are either in $E_R$ or $R$ or $W$. Within time 1 some process in $R$ or $W$ takes a step. If the first process not in $E_R$ taking a step started in $E_R$ or $R$, then it reaches $F$ and we are done; if the first process taking a step is in $W$, then it reaches $S$ since in $s$ no resource is held. Once a process $i$ is in $\{S, D, F\}$, then within two more time units process $i$ reaches either state $F$ or $P$, and we are done. ∎

18